\documentclass[a4paper]{article}
\usepackage[british]{babel}
\usepackage[utf8]{inputenc}
\usepackage[a4paper, total={6.5in, 9.5in}]{geometry}

\usepackage{amsmath,amssymb,amsthm}
\usepackage{dsfont}
\usepackage{tikz}
\usepackage{tikz-cd}
\usepackage[style=alphabetic,sorting=anyt,sortcites=true,backend=biber,maxbibnames=5,maxcitenames=5]{biblatex}
\usepackage{csquotes}
\usepackage{hyperref}
\usepackage{cleveref}
\usepackage{float}
\usepackage[table]{colortbl}% http://ctan.org/pkg/xcolor
\usepackage{longtable}
\usepackage{enumerate}

\title{Reidemeister zeta functions of low-dimensional almost-crystallographic groups are rational}
\author{Karel Dekimpe, Sam Tertooy, Iris Van den Bussche\thanks{Research supported  by long term structural funding -- Methusalem grant of the Flemish Government.}\\
KU Leuven Campus Kulak Kortrijk, E.~Sabbelaan 53, 8500 Kortrijk}
\date{\today}

\theoremstyle{plain}
\newtheorem{theorem}{Theorem}[section]
\newtheorem{lemma}[theorem]{Lemma}
\newtheorem{corollary}[theorem]{Corollary}

\theoremstyle{definition}
\newtheorem{defn}[theorem]{Definition}

\theoremstyle{remark}

\newtheorem{remark}[theorem]{Remark}

\DeclareMathOperator{\End}{End}
\DeclareMathOperator{\Aut}{Aut}

\DeclareMathOperator{\GL}{GL}

\DeclareMathOperator{\Aff}{Aff}

\DeclareMathOperator{\lcm}{lcm}
\newcommand{\I}{\mathds{1}}

\newcommand{\NN}{{\mathbb{N}}}
\newcommand{\ZZ}{{\mathbb{Z}}}

\newcommand{\RR}{{\mathbb{R}}}
\newcommand{\CC}{{\mathbb{C}}}
\newcommand{\FF}{{\mathbb{F}}}

\DeclareMathOperator{\Spec}{Spec}

\newcommand{\cc}[1]{{\bar{#1}}}
\newcommand{\lie}[1]{{\mathfrak {#1}}}
\begin{document}
	
	\maketitle
	
		\begin{center}
	This is an Accepted Manuscript of an article published by Taylor \& Francis in Communications in Algebra on 01 Mar 2018, available online:  \href{https://doi.org/10.1080/00927872.2018.1435792}{https://www.tandfonline.com/10.1080/00927872.2018.1435792}.
	\end{center}
	
\begin{abstract}
We prove that the Reidemeister zeta functions of automorphisms of crystallographic groups with diagonal holonomy \(\ZZ_2\) are rational. As a result, we obtain that Reidemeister zeta functions of automorphisms of almost-crystallographic groups up to dimension \(3\) are rational.
\end{abstract}

\section{The Reidemeister number and zeta function}
In this section we introduce basic notions concerning Reidemeister numbers. For a general reference on Reidemeister numbers and their connection with fixed point theory, we refer the reader to \cite{jian83-1}. In this paper, we use \(\NN\) to denote the set of positive integers and \(\NN_0\) to denote the set of non-negative integers. 
\begin{defn}
	Let \(G\) be a group and \(\varphi: G \to G\) an endomorphism. Define an equivalence relation \(\sim\) on \(G\) by
	\begin{equation*}
	\forall g,g' \in G: g \sim g' \iff \exists h \in G: g = hg'\varphi(h)^{-1}.
	\end{equation*}
	The equivalence classes are called \emph{Reidemeister classes} or \emph{twisted conjugacy classes}, and the number of equivalence classes is called the \emph{Reidemeister number} \(R(\varphi)\), which is therefore always a positive integer or infinity.
\end{defn}

\begin{defn}
	Let \(\Aut(G)\) be the automorphism group of a group \(G\). We define the \emph{Reidemeister spectrum} as
	\begin{equation*}
	\Spec_R(G) = \{R(\varphi) \mid \varphi \in \Aut(G)\}.
	\end{equation*}
	If \(\Spec_R(G) = \{\infty\}\) we say that \(G\) has the \emph{\(R_\infty\)-property}. 
\end{defn}

\begin{defn}
Let \(\varphi \in \End(G)\). The Reidemeister zeta function \(R_\varphi(z)\) of \(\varphi\) is defined as
\begin{equation*}
R_\varphi(z) = \exp \sum_{n=1}^{\infty} \frac{R(\varphi^n)}{n}z^n.
\end{equation*}
In order for this function to actually \emph{exist}, we require that \(R(\varphi^n) < \infty\) for all \(n \in \NN\), and that the series has a positive radius of convergence. 
\end{defn}

These Reidemeister zeta functions were introduced by A.~Fel'shtyn in  the late 80's. We refer the reader to \cite{fels00-2} for more details on the early results on these Reidemeister zeta functions. One of the central questions has always been whether or not such a Reidemeister zeta function is rational. A.~Fel'shtyn, together with R.~Hill, could already prove the rationality of this function in some special cases. E.g.\ for endomorphsims of finitely generated abelian groups and for finitely generated torsion-free nilpotent groups. In \cite{dd13-2} (see also \cite{fl15-1}), it was shown that the Reidemeister zeta function of any endomorphism (if this zeta function is defined) on a finitely generated torsion-free virtually nilpotent is rational. 

Up till now, the situation for virtually nilpotent groups which are not torsion-free has not really been studied and it is not sure that we should expect that also in this case all Reidemeister zeta functions will be rational. The aim of this paper is to take off with the study of Reidemeister zeta functions of automorphisms in this non torsion-free case. 
In the cases we treat, we do find that all Reidemeister zeta functions are indeed rational.

\section{Almost-crystallographic groups}
Let \(G\) be a connected and simply connected nilpotent Lie group. We define the group of affine transformations on \(G\) as the semi-direct product \(\Aff(G) = G \rtimes \Aut(G)\), where multiplication is defined by \((d_1,D_1)(d_2,D_2) = (d_1D_1(d_2),D_1D_2)\). 
Let \(C\) be a maximal compact subgroup of \(\Aut(G)\), then \(G \rtimes C\) is a subgroup of \(\Aff(G)\). A cocompact discrete subgroup \(\Gamma\) of \(G \rtimes C\) is called an \emph{almost-crystallographic group}, or an \emph{almost-Bieberbach group} if it is also torsion-free. If \(G = \RR^n\), then we call \(\Gamma\) a \emph{crystallographic group} or a \emph{Bieberbach group} respectively.

Crystallographic groups are well understood by the three Bieberbach theorems. We refer to  \cite{wolf77-1,char86-1,szcz12-1} for more information on and proofs of these theorems.  These theorems have been generalised to almost-crystallographic groups (see \cite{deki96-1} for more details).

\begin{theorem}[Generalised first Bieberbach theorem]
	Let \(\Gamma \leq \Aff(G)\) be an almost-crystallographic group. Then \(N = \Gamma \cap G\) is a uniform lattice of \(G\) and \(N\) is of finite index in \(\Gamma\). 
\end{theorem}

Moreover, \(N\) is the unique maximal nilpotent  and normal (hence characteristic) subgroup of \(\Gamma\). Hence
any almost-crystallographic group \(\Gamma\) fits in a short exact sequence
\begin{equation*}
1 \to N \to \Gamma \to F \to 1,
\end{equation*}
with \(F\) a finite group called the \emph{holonomy group} of \(\Gamma\).  Note that \(F=p(\Gamma)\), where \(p:\Aff(G)=G\rtimes \Aut(G)\to \Aut(G)\) is the natural projection onto the second factor. So, we can (and will) view \(F\) as being a subgroup of \(\Aut(G)\).
\begin{theorem}[Generalised second Bieberbach theorem]
	\label{thm:gensecbieb}
	Let \(\Gamma \leq \Aff(G)\) be an almost-crystallographic group and \(\varphi \in \Aut(\Gamma)\). Then there exists some \((d,D) \in \Aff(G)\) such that \(\varphi(\gamma) = (d,D)\gamma(d,D)^{-1}\) for all \(\gamma \in \Gamma\).
\end{theorem}

Let \(\Gamma\leq \Aff(\RR)\) be a crystallographic group, then \(\Gamma\cap \RR^n\) is a uniform lattice of \(\RR^n\) and hence is isomorphic to \(\ZZ^n\). After an affine conjugation we may suppose that \(\Gamma\cap \RR^n=\ZZ^n\) (and not just isomorphic to \(\ZZ^n\)). In this case \(\Gamma \leq \RR^n \rtimes \GL_n(\ZZ)\) (where we view $\GL_n(\ZZ)$ as a subgroup of \(\GL_n(\RR)\)). It follows that \(F\) is a subgroup of \(\GL_n(\ZZ)\), and therefore, the matrix \(D\) from \cref{thm:gensecbieb} will belong to the normaliser \(N_{\GL_n(\ZZ)}(F)\). We will call the holonomy group \(F\) \emph{diagonal} if it consists entirely of diagonal matrices.

\medskip

For any morphism of Lie groups \(\Phi: G \to H\), we will denote the induced morphism on the corresponding Lie algebras by means of a star-index: \(\Phi_*: \lie{g} \to \lie{h}\). Note that when \(G=H=\RR^n\), then \(\Phi\) can be seen as an element of \(\GL_n(\RR)\) and in that case \(\Phi_\ast=\Phi\). 

\medskip

We will need the following criterion:
\begin{theorem}[see \cite{dp11-1}]
	\label{thm:det1-AD}
	Let \(\Gamma\) be an almost-crystallographic group with holonomy group \(F\). Let \(\varphi: \Gamma \to \Gamma\) be an automorphism such that \(\varphi(\gamma) = (d,D)\gamma(d,D)^{-1}\) for every \(\gamma \in \Gamma\), with \((d,D) \in \Aff(G)\). Then
	\begin{equation*}
	R(\varphi) = \infty \iff \exists A \in F \mbox{ such that } \det(\I - A_*D_*) = 0.
	\end{equation*}
\end{theorem}

\section{Existence and rationality of Reidemeister zeta functions}
The goal of this section is to determine which almost-crystallographic groups admit Reidemeister zeta functions, and for which of those groups we already know the rationality of the zeta functions.

First of all, it is obvious that a group admitting the \(R_\infty\)-property does not admit any Reidemeister zeta function for automorphisms. However, other criteria do exist, as demonstrated by the following theorem.
\begin{theorem}
	Let \(\Gamma\) be an almost-crystallographic group with a characteristic subgroup \(H \cong \ZZ\). Then no Reidemeister zeta function of automorphisms of \(\Gamma\) exists.
\end{theorem}
\begin{proof}
	Let \(H = \langle x \rangle\). As \(H\) is normal and abelian, we must have that \(H\leq N =\Gamma\cap G\). (Recall that \(N\) is the unique maximal nilpotent and normal subgroup of \(\Gamma\)).  Since \(N\) is nilpotent and \(H\) is normal in \(N\), we must have that the intersection \(H\cap Z(N)\) of \(H\) with the centre of \(N\) is non-trivial. So, there exists some \(k>0\) such that \(x^k\in Z(N)\). In fact, as \(N\) is torsion-free, \(N/Z(N)\) is torsion-free as well and so we have that \(x\in Z(N)\), hence \(H\leq Z(N)\). Let \(\varphi\) be an automorphism of \(\Gamma\) given by \(\varphi(\gamma) = (d,D)\gamma(d,D)^{-1}\). As $x\in Z(N)$, it then follows that $\varphi(x)=D(x)$. Either \(\varphi(x) = D(x)= x\) or \(\varphi(x) = D(x)= x^{-1}\). In any case we have that $D^2(x)=x$. It then follows that there exists a non-zero element $X\in \lie{g}$ in the Lie algebra of $G$ with $D^2_\ast(X)=X$  and therefore \(\det(\I-D^2_*) = 0\). So certainly \(R(\varphi^2) = \infty\) and we can conclude that the Reidemeister zeta function \(R_\varphi(z)\) does not exist.
\end{proof}
A particular application of this theorem becomes clear when we recall the following property of three-dimensional almost-crystallographic groups (that are not crystallographic). 
\begin{lemma}[see {\cite[Lemma 4.3]{dp11-1}}]
	Let \(\Gamma\) be an almost-crystallographic group with \(N\) isomorphic to the Heisenberg group. Then \(\Gamma \cap Z(N)\) is a characteristic subgroup of \(\Gamma\) that is isomorphic to \(\ZZ\). 
\end{lemma}
Therefore, we may limit ourselves to crystallographic groups, and we may of course exclude the crystallographic groups with centre isomorphic to \(\ZZ\). Of the few crystallographic groups remaining that do not have the \(R_\infty\)-property, we may still exclude some using the next lemma.
\begin{lemma}
	\label{lem:crystgroupsnozeta}
	Let \(\Gamma\) be an \(n\)-dimensional crystallographic group with holonomy group \(F \leq \GL_n(\ZZ)\). If \(|N_{\GL_n(\ZZ)}(F)| < \infty\), then no Reidemeister zeta functions of automorphisms of \(\Gamma\) exist.
\end{lemma}

\begin{proof}
	Suppose that \(|N_{\GL_n(\ZZ)}(F)| < \infty\). Any automorphism \(\varphi\) can be written as \(\varphi(\gamma) = (d,D)\gamma (d,D)^{-1}\) for all \(\gamma \in \Gamma\), with \(d \in \RR^n\) and \(D \in N_{\GL_n(\ZZ)}(F)\); and \(R(\varphi) < \infty\) if and only if \(|\det(\I_n-AD)| \neq 0\) for all \(A \in F\). Since the normaliser of the holonomy group is finite, there exists some power \(k\) such that \(D^k = \I_n\). Then \(|\det(\I_n-D^k)| = 0\), hence \(R(\varphi^k) = \infty\) and thus the Reidemeister zeta function of \(\varphi\) does not exist.
\end{proof}

Finally, we may also skip (almost-)Bieberbach groups, since rationality is already known there.

\begin{theorem}[see {\cite[Corollary 4.7]{dd13-2}}]
	Reidemeister zeta functions of almost-Bieberbach groups are rational.
\end{theorem}

In fact, in \cite{dd13-2} the rationality is proved for the Nielsen zeta function, but it is known that when the Reidemeister zeta function is defined for an automorphism, then it coincides with a Nielsen function of a map on the corresponding infra-nilmanifold (see \cite[Proposition 3.2]{fl15-1}).

\medskip

Thus the only almost-crystallographic groups that are left to check, are \(\langle \ZZ^2, (0,-\I_2)\rangle\) and \(\langle \ZZ^3, (0,-\I_3)\rangle\), i.e. the subgroup of \(\Aff(\RR^n)\) (\(n= 2,3\)) generated by \((0,-\I_n)\) and the elements \((x,\I_n)\) with \(x \in \ZZ^n\). Rather than limit ourselves to these two groups, we will prove the rationality of the Reidemeister zeta functions of all crystallographic groups with diagonal holonomy \(\ZZ_2\).

\section{The crystallographic groups \(\langle \ZZ^n, (0,-\I_n)\rangle\)}
\label{sec:Zn0-1}
The first step in calculating a Reidemeister zeta function is, of course, calculating the Reidemeister numbers \(R(\varphi^k)\). The following lemma can be found in \cite[Proposition 5.10]{dkt17-2}.
\begin{lemma}\label{lemma.4.1}
	Let \(\Gamma = \langle \ZZ^n, (0,-\I_n)\rangle\) and \(\varphi \in \Aut(\Gamma)\). Then there exists an affine map \((\frac{d}{2}, D)\in \Aff(\RR^n)\) with \(d\in \ZZ^n\) such that \(\varphi(\gamma)\) is given by \(\varphi(\gamma) = (d,D)\gamma(d,D)^{-1}\). If \(R(\varphi) < \infty\), then
	\begin{equation}
	\label{eq:Rformula}
	R(\varphi) = \frac{|\det(\I_n - D)|+|\det(\I_n + D)|}{2} + O(\I_n - D,d),
	\end{equation}
	where for any \(A\in \ZZ^{n \times n}\) and \(a \in \ZZ^n\), we use \(O(A,a)\) to denote the number of solutions \(x\) over \(\ZZ_2\) of the linear system \(\bar{A}x = \bar{a}\). Here we use \(\bar{A}\)  $($resp.~$\bar{a})$ to denote the natural projection of \(A\) $($resp.~\(a)\) to \(\ZZ_2^{n \times n}\) $($resp.~$\ZZ_2^n)$.
\end{lemma}
Let \(\lambda_1, \dots, \lambda_n\) be the eigenvalues of \(D\). We may then rewrite \eqref{eq:Rformula} as
\begin{align*}
R(\varphi) &= \frac{\left|\prod_{i=1}^{n}(1-\lambda_i)\right|+\left|\prod_{i=1}^{n}(1+\lambda_i)\right|}{2} + O(\I_n - D,d)\\
&=\frac{\prod_{i=1}^{n}|1-\lambda_i|+\prod_{i=1}^{n}|1+\lambda_i|}{2} + O(\I_n - D,d).
\end{align*}
Similarly, for any \(k \in \NN\) we have
\begin{equation*}
R(\varphi^k) = \frac{\prod_{i=1}^{n}|1-\lambda_i^k|+\prod_{i=1}^{n}|1+\lambda_i^k|}{2} + O\left(\I_n-D^k,\left[\sum_{i=0}^{k-1}D^i\right] d\right).
\end{equation*}
We will deal with both terms separately in the following subsections.

\subsection{The first term}
We will need the following lemma to simplify the first term.
\begin{lemma}
	\label{lem:firsttermkthpowers}
	Let \(\lambda_1, \dots, \lambda_n\) be the eigenvalues of some matrix \(D \in \GL_n(\ZZ)\). Then there exist non-negative  integers \(a, b \in \NN_0\) and complex numbers \(\mu_1, \dots, \mu_{a},\nu_{1}, \dots, \nu_b\) such that 
	\begin{equation*}
	\frac{\prod_{i=1}^{n}|1-\lambda_i^k|+\prod_{i=1}^{n}|1+\lambda_i^k|}{2} = \mu_1^k + \mu_2^k + \cdots + \mu_{a}^k - \nu_{1}^k - \cdots - \nu_b^k
	\end{equation*}
	for each \(k \in \NN\).
\end{lemma}
\begin{proof}
	The eigenvalues \(\lambda_i\) can either be real or complex (and non-real), and if some eigenvalue \(\lambda_i\) is complex, then its complex conjugate \(\cc{\lambda}_i\) will also be an eigenvalue.
	
	Let us first consider a real eigenvalue \(\lambda_i\). We distinguish three cases:
	\begin{enumerate}[(1)]
		\item \(|\lambda_i| < 1\). Then \(|1-\lambda_i^k| = 1^k-\lambda_i^k\) and \(|1+\lambda_i^k| = 1^k+\lambda_i^k\),
		\item \(\lambda_i > 1\). Then \(|1-\lambda_i^k| = -1^k+\lambda_i^k\) and \(|1+\lambda_i^k| = 1^k+\lambda_i^k\),
		\item \(\lambda_i < -1\). Then \(|1-\lambda_i^k| = -(-1)^k+(-\lambda_i)^k\) and \(|1+\lambda_i^k| = (-1)^k+(-\lambda_i)^k\).
	\end{enumerate}
	So for real \(\lambda_i\), both \(|1-\lambda_i^k|\) and \(|1+\lambda_i^k|\) can be written as the sum and/or difference of \(k\)-th powers.
	
	Next, let us consider a complex eigenvalue \(\lambda_i\). Then its complex conjugate \(\cc{\lambda}_i\) is also an eigenvalue, and taking their respective factors together we get:
	\begin{equation*}
	|1 \pm \lambda_i^k||1 \pm \cc{\lambda}_i^k| = |1\pm \lambda_i^k|^2 = 1^k  \pm \lambda_i^k  \pm \cc{\lambda}_i^k  + (|\lambda_i|^2)^k.
	\end{equation*}
	The first equality tells us that this product is real and moreover positive, the second equality tells us that this can be written as the sum and/or difference of \(k\)-th powers.

	Combining the real and complex cases, we can expand both products \(\prod_{i=1}^{n}|1-\lambda_i^k|\) and \(\prod_{i=1}^{n}|1+\lambda_i^k|\) and obtain a sum of terms of the form $\pm\lambda_{i_1}^k\lambda_{i_2}^k \cdots \lambda_{i_p}^k= \pm (\lambda_{i_1}\lambda_{i_2} \cdots \lambda_{i_p})^k$ (where $p$ varies between 0 and $n$). Note that all of these terms are, up to sign,   \(k\)-th powers of terms which themselves do not depend on $k$. These two products will have exactly the same terms, though the sign of said terms may differ. If two matching terms have the same sign, their sum will have a factor \(2\) that cancels out with the \(2\) in the denominator; and if two matching terms have the opposite sign, they cancel out each other. So the entire term is indeed a sum and/or difference of \(k\)-th powers of fixed terms (not depending on $k$).
\end{proof}
With this lemma proven, it is now easy to show the rationality of the first term.
\begin{lemma}
	\label{lem:part1}
	Let \(\lambda_1, \dots, \lambda_n\) be the eigenvalues of some matrix \(D \in \GL_n(\ZZ)\). The function
	\begin{equation*}
	\exp \sum_{k=1}^{\infty}
	\frac{\prod_{i=1}^{n}|1-\lambda_i^k|+\prod_{i=1}^{n}|1+\lambda_i^k|}{2}\frac{z^k}{k}
	\end{equation*}
	is a rational function.
\end{lemma}
\begin{proof}
	We invoke the previous lemma to obtain
	\begin{align*}
	\exp \sum_{k=1}^{\infty}
\frac{\prod_{i=1}^{n}|1-\lambda_i^k|+\prod_{i=1}^{n}|1+\lambda_i^k|}{2}\frac{z^k}{k}
	&= \exp \sum_{k=1}^{\infty}\frac{z^k}{k} \left(\sum_{i=1}^{a} \mu_i^k - \sum_{i=1}^b \nu_i ^k\right)\\
	&= \exp \left(\sum_{i=1}^{a}\sum_{k=1}^{\infty} \frac{\mu_i^k}{k}z^k - \sum_{i=1}^{b}\sum_{k=1}^{\infty} \frac{\nu_i^k}{k}z^k  \right)\\
	&= \exp \left(- \sum_{i=1}^{a}\log(1-\mu_i z) +  \sum_{i=1}^{b}\log(1-\nu_i z) \right)\\
	&= \frac{ \prod_{i=1}^{b}(1-\nu_i z)}{\prod_{i=1}^{a}(1-\mu_i z)},
	\end{align*}
	which is a rational function.
\end{proof}

\subsection{The second term}
The second term is far less straightforward. We first introduce a particular family of sequences.
\begin{defn}
	We define the sequence \(a^i = (a^i_k)_{k \in \NN}\) by
	\begin{equation*}
	a^i_k =  \begin{cases}
	i & \text{ if }k \equiv 0 \bmod i,\\
	0 & \text{otherwise.}
	\end{cases}
	\end{equation*}
\end{defn}
The vast majority of this subsection will be devoted to proving the following theorem.
\begin{theorem}
	\label{thm:Oissumsequences}
	Let \(D \in \GL_n(\ZZ), d \in \ZZ^n\). Then there exist \(l \in \NN_0\) and \(c_1, \dots, c_l \in \NN_0\) such that
	\begin{equation*}
	 O\left(\I_n-D^k,\left[\sum_{i=0}^{k-1}D^i\right] d\right)  = c_1 a^1_k + c_2 a^2_k + \cdots + c_l a^l_k
	\end{equation*}
	for all \(k \in \NN \).
\end{theorem}

The actual proof of this theorem will require some preparation. 

\medskip

Before we really start with the  proof of this theorem, let us note that we do not need full information on the pair $(d,D)$, but we only need to know their natural projections modulo 2, namely the pair $(\bar{d}, \bar{D})$ as mentioned in the formulation of \cref{lemma.4.1}. To avoid having to write a bar above $d$ and $D$ each time we will assume from now onwards that $D\in \GL_n(\ZZ_2)$ and $d\in \ZZ_2^n$.

\medskip

We will apply a change of base such that \(D\) has a more suitable form to work with. With that in mind, we first need the following matrix decomposition.

\begin{lemma}
	\label{lem:NXXDB}
	Let \(N\) be a nilpotent, upper-triangular \(k \times k\)-matrix and \(D\) an invertible \(l \times l\)-matrix over a field \(\FF\). For any \(k \times l\)-matrix \(B\), there exists a (unique) \(k \times l\)-matrix \(X\) such that
	\begin{equation}
	\label{eq:NXXDB}
	NX + XD = B.
	\end{equation}
\end{lemma}
\begin{proof}
	We prove this by induction on \(k\). If \(k = 1\), then \(N = 0\) and \(X = BD^{-1}\). Now let \(k \geq 2\) and suppose that the lemma holds for smaller values of \(k\). Then \(N\), \(X\) and \(B\) can be seen as block matrices of the forms
	\begin{equation*}
	N = \left(\begin{array}{ccc|c}
	&&&\\
	&N_1 & & N_2\\
	&&&\\
	\hline
	0 & \cdots & 0 & 0
	\end{array}\right),\quad X =  \left(\begin{array}{ccc}
	&&\\
	&X_1&\\
	&&\\
	\hline
	 & X_2 &
	\end{array}\right),\quad B =  \left(\begin{array}{ccc}
	&&\\
	&B_1&\\
	&&\\
	\hline
	& B_2 &
	\end{array}\right),
	\end{equation*}
	where \(N_1\) is a nilpotent, upper-triangular \((k-1) \times (k-1)\)-matrix, \(N_2\) is a \((k-1)\times 1\)-matrix, \(X_1\) and \(B_1\) are \((k-1)\times l\)-matrices and \(X_2\) and \(B_2\) are \(1\times l\)-matrices. We can then split up \eqref{eq:NXXDB} in the system of equations
	\begin{equation*}
	\begin{cases}
	N_1 X_1 + N_2 X_2 + X_1 D = B_1,\\
	X_2 D = B_2.
	\end{cases}
	\end{equation*}
	The second equation gives us \(X_2 = B_2D^{-1}\), and substituting this into the first equation gives
	\begin{equation*}
	N_1X_1 + X_1 D = B_1 - N_2B_2D^{-1}.
	\end{equation*}
	By applying the induction hypothesis, we get a solution \(X_1\). Together with \(X_2\) we have the full solution \(X\) of \eqref{eq:NXXDB}.
\end{proof}
This decomposition allows us to put \(D\) in the required form.
\begin{lemma}
	\label{lem:Dblocks}
	Let \(D\) be an \(n \times n\)-matrix over a field \(\FF\), then there exists an invertible matrix \(P\) such that
	\begin{equation*}
	PDP^{-1} = \left(\begin{array}{cc}
	D_1 & 0\\0 & D_2
	\end{array}\right),
	\end{equation*}
	where \(D_1\) is a unipotent, upper-triangular matrix, and \(D_2\) does not have eigenvalue \(1\) (and hence \(\I-D_2\) is invertible).
\end{lemma}
\begin{proof}
	Consider the linear map \begin{equation*}
	f: \FF^n \to \FF^n: \vec{x} \mapsto D\vec{x}.
	\end{equation*}
	It suffices to show that there exists a basis such that \(f\) has the required form with respect to this basis. Suppose that \(D\) has eigenvalue \(1\), then take an eigenvector corresponding to this eigenvalue and extend to a basis. With respect to this basis, we have
	\begin{equation*}
	D \sim \left(\begin{array}{c|ccc}
	1 & * & \cdots & *\\
	\hline
	0 & & & \\
	0 & & D'&\\
	0 & & &
	\end{array}\right).
	\end{equation*}
	We can then interpret \(D'\) as a linear map \(\FF^{n-1} \to \FF^{n-1}\) and proceed by induction to obtain
	\begin{equation*}
		D \sim \left(\begin{array}{c|c}
		D_1 & B\\
		\hline
		0 & D_2
		\end{array}\right),
	\end{equation*}
	with \(D_1\) a unipotent upper-triangular \(k \times k\)-matrix and \(D_2\) an \(l \times l\)-matrix with no eigenvalue \(1\).  Hence \(D_1-\I_k\) is a nilpotent upper-triangular \(k \times k\)-matrix and \(\I_l - D_2\) is an invertible \(l \times l\)-matrix. By \cref{lem:NXXDB} there exists a \(k \times l\)-matrix \(X\) such that
	\begin{equation*}
	(D_1 - \I_k) X + X(\I_l-D_2) = B,
	\end{equation*}
	which in turn gives
	\begin{equation*}
	\begin{pmatrix}
	\I_k & X\\
	0 & \I_l
	\end{pmatrix}\begin{pmatrix}
	D_1 & B\\
	0 & D_2
	\end{pmatrix}\begin{pmatrix}
	\I_k & X\\
	0 & \I_l
	\end{pmatrix}^{-1}
	= \begin{pmatrix}
	D_1 & 0\\
	0 & D_2
	\end{pmatrix},
	\end{equation*}
	as required.
\end{proof}
From this point onwards, we will work with \(\FF = \ZZ_2\). To any pair \((d,D)\), with $d\in \ZZ_2^n$ and $D\in \GL_n(\ZZ_2)$, we associate the sets \(V_k\) and \(W_k\) defined as
\begin{align*}
V_k &= \left\{x \in \ZZ_2^n  \;\middle|\; (\I - D^k)x = \left[\sum_{i=0}^{k-1}D^i\right] d \right\},\\
W_k &= \left\{ x \in V_k  \;|\; x \notin V_l \quad \forall l \in \{1,2, \dots, k-1\}\right\}.
\end{align*}
Let \(v_k = |V_k| = O\left(\I_n-D^k,\left[\sum_{i=0}^{k-1}D^i\right] d\right)\) and \(w_k = |W_k|\). The $W_k$ are disjoint sets and their union is all of $\ZZ_2^n$. Hence, it is obvious that only for a finite number of values of $k$ we have that $w_k\neq 0$,  since their sum equals $2^n$. To prove \cref{thm:Oissumsequences}, we need to determine what the sequence  \(v = (v_k)_{k \in \NN}\) is. As we have split up \(D\) in a unipotent block \(D_1\) and a block with no eigenvalue \(1\), \(D_2\), we will first restrict to these two blocks.

\subsubsection{If \(D\) has no eigenvalue \(1\).}
Let us first assume that \(D\) does not have eigenvalue \(1\), and therefore \(\I-D\) is invertible. Then there exists some \(d_0\) such that \((\I-D)d_0 = d\), and hence we can state
\begin{equation*}
\left[\sum_{i=0}^{k-1}D^i\right] d =   \left[\sum_{i=0}^{k-1}D^i\right] (\I - D)d_0 = (\I-D^k)d_0,
\end{equation*}
so we are actually searching for solutions of the linear system given by
\begin{equation*}
(\I - D^k)(x-d_0) = 0.
\end{equation*}
The ``shift'' by \(d_0\) has no effect on the number of solutions of this system, so we may assume without loss of generality that
\begin{equation*}
V_k = \left\{x \in \ZZ_2^n  \;|\; (\I - D^k)x = 0\right\}.
\end{equation*}
We will now formulate and prove some properties of these sets \(V_k\) and \(W_k\).
\begin{lemma}
	\label{lem:vkvl=>vl-k}
	Let \(k < l\). If \(x \in V_k \cap V_l\), then \(x \in V_{l-k}\).
\end{lemma}
\begin{proof}
	Let \(x \in V_k \cap V_l\). Then
	\begin{equation*}
	0 = (\I-D^l)x = (\I-D^k +D^k - D^l)x = (\I - D^k)x + D^k(\I - D^{l-k})x = D^k(\I - D^{l-k})x,
	\end{equation*}
	and because \(D\) is invertible we are left with \((\I - D^{l-k})x = 0\), hence \(x \in V_{l-k}\).
\end{proof}
\begin{corollary}
	Let \(k < l\). If \(x \in V_k \cap V_l\), then\begin{enumerate}[(1)]
		\item \( x \in V_{l \bmod k}\),
		\item \( x \in V_{\gcd(k,l)}\).
	\end{enumerate}
\end{corollary}
\begin{proof}
	The first property follows by repeatedly applying \cref{lem:vkvl=>vl-k}. The second property follows by repeatedly applying the first property.
\end{proof}
On the other hand, we also have
\begin{lemma}
	If \(x \in V_k\), then \(x \in V_{kl}\) for all \(l \in \NN\).
\end{lemma}
\begin{proof}All we have to do is split \((\I-D^{kl})\) in suitable factors:
	\begin{equation*}
	(\I-D^{kl})x = (\I + D^k + \dots + D^{(l-1)k})(\I - D^k)x = 0,
	\end{equation*}
	because \((\I - D^k)x = 0\).
\end{proof}

In conclusion, we can state that \(V_k\) is exactly the disjoint union
\begin{equation*}
V_k = \bigsqcup_{d | k} W_d.
\end{equation*}
and hence we get
\begin{equation*}
v_k = \sum_{d|k}w_d = \sum_{d|k} \frac{w_d}{d}(a^d)_k,
\end{equation*}
where we used that the \(k\)-th element in the sequence \(a^d\) is \(d\), since \(k\) is a multiple of \(d\). 
Now, since $(a^d)_k=0$ when $d$
 does not divide $k$ we have that 
\[ v_k = \sum_{d=1}^k \frac{w_d}{d}(a^d)_k= \sum_{d=1}^\infty \frac{w_d}{d}(a^d)_k.\]

 So we indeed seem to have a sum of sequences \(a^d\), but we still require the coefficients of this sum to be integers.

\begin{lemma}
	\label{lem:kdivideswk}
	 \(k\) divides \(w_k\) for any \(k \in \NN\).
\end{lemma}
\begin{proof}
We define an action of \(\ZZ\) on \(W_k\) by
\begin{equation*}
\ZZ \times W_k \to W_k: (z,x) \mapsto z \cdot x = D^z x.
\end{equation*}
First, we verify that this action is well-defined. If \(x \in W_k\), then 
\begin{equation*}
(\I - D^k)D^zx = D^z(\I-D^k)x = 0,
\end{equation*}
hence \(D^zx \in V_k\). On the other hand, if for some \(l < k\) we were to have that \(D^zx  \in V_l\), then 
\begin{equation*}
0 = (\I-D^l)D^zx = D^z(\I-D^l)x.
\end{equation*}
Because \(D\) is invertible, this would mean that \((\I-D^l)x = 0\), or in other words \(x \in V_l\). This is a contradiction since \(x \in W_k\). In fact, \(k\ZZ\) acts trivially on \(W_k\) since 
\begin{equation*}
(\I-D^k)x = 0 \iff D^k x = x,
\end{equation*}
so we can redefine the original action as an action of \(\ZZ_k\) on \(W_k\), which is a free action. Indeed, suppose that for some \(x \in W_k\) we have that \(D^lx = x\), where \(l\) is not a multiple of \(k\). Then \(x \in V_l\) and therefore \(x \in V_{l \bmod k}\). This obviously contradicts that \(x \in W_k\).

By the orbit-stabiliser theorem, we can now partition \(W_k\) into finitely many orbits of length \(k\), and thus \(k\) divides \(w_k\). 
\end{proof}
Putting everything together now, we can conclude that the sequence \(v = (v_k)_{k \in \NN}\) equals
\begin{equation*}
v = \sum_{k=1}^\infty \frac{w_k}{k} a^k,
\end{equation*}
which has integer coefficients since \(k|w_k\). Recall that this is  actually a finite sum, since only finitely many of the \(w_k\) are non-zero.

\subsubsection{If \(D\) is unipotent upper-triangular.}
For the case where \(D\) is unipotent upper-triangular, we will have very similar results as the previous case. The main difference here will be that we will end up working mainly with powers of \(2\) as opposed to arbitrary \(k\). Because we are working over \(\ZZ_2\), we have the following two statements:
\begin{remark}
	\label{rem:suminvertible}
	If \(m\) is an odd positive integer, then for any integers \(k_1, k_2, \dots, k_m\), we have that \(D^{k_1} + D^{k_2} + \cdots + D^{k_m}\) is unipotent upper-triangular (and hence invertible).
\end{remark}

\begin{remark}
	If \(D\) is an \(n \times n\)-matrix, then \(D^{2^{n-1}} = \I_n\), since \(\I_n - D^{2^{n-1}} = (\I_n - D)^{2^{n-1}} = 0\). This means that \(V_{2^n} = \ZZ_2^n\).
\end{remark}

The next lemma makes clear why we only really need to care about powers of \(2\).
\begin{lemma}
	Decompose \(k\) as \(k = 2^rm\) with \(m\) odd. Then \(V_k = V_{2^r}\).
\end{lemma} 
\begin{proof}
	Let \(M = \I + D^{2^r} + D^{2 \cdot 2^r} + \cdots + D^{(m-1)2^r}\), which is invertible (see \cref{rem:suminvertible}). Then 
	\begin{equation*}
	\I - D^k =  \I - D^{2^rm} = (\I + D^{2^r} + D^{2 \cdot 2^r} + \cdots + D^{(m-1)2^r})(\I - D^{2^r}) = M(\I - D^{2^r}),
	\end{equation*}
	and
	\begin{equation*}
	\left[\sum_{i=0}^{k-1}D^i \right] d = \left[\sum_{i=0}^{2^rm-1}D^i \right]d = M \left[\sum_{i=0}^{2^r-1}D^i\right] d.
	\end{equation*}
	We then obtain
	\begin{align*}
	(\I - D^k)x = \left[\sum_{i=0}^{k-1}D^i \right] d &\iff M (\I - D^{2^r})x = M \left[\sum_{i=0}^{2^r-1}D^i\right] d\\
	& \iff (\I - D^{2^r})x =  \left[\sum_{i=0}^{2^r-1}D^i\right] d,
	\end{align*}
	and therefore \(V_k = V_{2^r}\).
\end{proof}
We conclude that \(w_k = 0\) if \(k \) is not a power of \(2\). Now let  \(r_0\) be the smallest power of \(2\) such that \(w_{2^{r_0}} \neq 0\), then there exists some \(d_0\) such that 
\begin{equation*}
(\I - D^{2^{r_0}})d_0 = \left[\sum_{i=0}^{2^{r_0}-1} D^i\right] d.
\end{equation*}
Similarly to the other case, the number of elements $v_{2^{r_0}}=w_{2^{r_0}}$ in $V_{2^{r_0}}=W_{2^{r_0}}$ is equal to the number of solutions of the system
\begin{equation*}
(\I-D^{2^{r_0}})(x-d_0) = 0.
\end{equation*}

\begin{lemma}
	\label{lem:2r0dividesw2r0}
	\(2^{r_0}\) divides \(w_{2^{r_0}}\).
\end{lemma}
\begin{proof}
	As we are working over \(\ZZ_2\), we have that \(\I-D^{2^{r_0}} = (\I-D)^{2^{r_0}}\). Since \(D\) is unipotent upper-triangular, \(\I-D\) is nilpotent upper-triangular, and hence taking the \(2^{r_0}\)-th power gives a matrix where the bottom \(r_0\) rows are zero. Thus \(w_{r_0} = |\ker(\I-D^{2^{r_0}})|\) is a multiple of \(2^{r_0}\).
\end{proof}
We have already shown that if \(r = r_0\), we may work with the linear system \((\I-D^{2^{r}})(x-d_0) = 0\). This is, however, rather useless if we do not have this for every \(r\). For \(r > r_0\) we have
\begin{align*}
(\I-D^{2^r})x &= \left[\sum_{i=0}^{2^{r}-1} D^i\right] d\\
&= (\I + D^{2^{r_0}} +  D^{2 \cdot 2^{r_0}} + \cdots +  D^{(2^{r-r_0}-1)2^{r_0}})\left[\sum_{i=0}^{2^{r_0}-1} D^i\right] d\\
&= (\I + D^{2^{r_0}} +  D^{2 \cdot 2^{r_0}} + \cdots +  D^{(2^{r-r_0}-1)2^{r_0}})(\I-D^{2^{r_0}})d_0\\
&=(\I-D^{2^r})d_0.
\end{align*}
So indeed we end up with the linear system
\begin{equation*}
(\I-D^{2^{r}})(x-d_0) = 0,
\end{equation*}
and again we may assume without loss of generality that \(d_0 = 0\). Now that we have this system for all \(r \geq r_0\), we also want to generalise \cref{lem:2r0dividesw2r0} to all \(r > r_0\). 

\begin{lemma}
	\(2^{r}\) divides \(w_{2^{r}}\) for all \(r > r_0\).
\end{lemma}
\begin{proof}
	Analogously to \cref{lem:kdivideswk}, we have an action of \(\ZZ_{2^r}\) on \(W_{2^r}\). Suppose this action is not free. As subgroups of \(\ZZ_{2^r}\) are generated by divisors of \(2^r\), there then exist some \(x \in W_{2^r} \) and some \(r' < r\) such that \(D^{2^{r'}}x = x\), which contradicts that \(x \in W_{2^r}\). So \(2^r \) divides \(w_{2^r}\).
\end{proof}

The following steps are identical to the case where \(D\) has no eigenvalue \(1\), hence we leave these to the reader and we can conclude that also in this case

\[v = \sum_{k=1}^\infty \frac{w_k}{k} a^k ,\]
where $\frac{w_k}{k}$ is an integer and the sum is in fact finite.

\subsubsection{The general case}
We now have all the necessary tools to prove \cref{thm:Oissumsequences}.
\begin{proof}[Proof of \cref{thm:Oissumsequences}]
	At the start of this section we proved that, after a change of basis, \(D\) is a block matrix of the form
	\begin{equation*}
	D = \begin{pmatrix}
	D_1 & 0\\
	0 & D_2
	\end{pmatrix},
	\end{equation*}
	such that \(D_1\) is unipotent upper-triangular and \(D_2\) has no eigenvalue \(1\). We may split the vector \(d\) in two pieces \(d_1\) and \(d_2\) matching the sizes of \(D_1\) and \(D_2\) respectively. So for any \(k\), we have two linear systems of equations given by
	\begin{equation*}
	\begin{cases}
	(\I - D_1^k)x_1 = \left[\sum_{i=0}^{k-1}D_1^i\right] d_1, \\
	(\I - D_2^k)x_2 = \left[\sum_{i=0}^{k-1}D_2^i\right] d_2 .
	\end{cases}
	\end{equation*}
	The total number of solutions \(x\) is of course the number of pairs \((x_1,x_2)\). In the previous subsections we have shown that both ``subsystems'' give sequences
	\begin{align*}
	v = (v_1,v_2,v_3,\dots), \\
	v' = (v'_1,v'_2,v'_3,\dots),
	\end{align*}
	that are linear combinations of the sequences \(a^i\), say
	$v=\sum_{k=1}^n c_k a^k$ and $v'= \sum_{l=1}^m c'_l a^l$. To solve the linear system as a whole, we are actually looking for the sequence \(v\cdot v'\) given by the component-wise multiplication of \(v\) and \(v'\):
	\begin{equation*}
	v \cdot v' = (v_1v'_1,v_2v'_2,v_3v'_3,\dots).
	\end{equation*}
	Using that \(a^k \cdot a^l = \gcd(k,l)a^{\lcm(k,l)}\), we get
	\begin{equation*}
	v \cdot v' = \left(\sum_{k=1}^n c_k a^k\right)\left(\sum_{l=1}^m c'_l a^l\right) = \sum_{k,l} c_kc'_l \gcd(k,l)a^{\lcm(k,l)},
	\end{equation*}
	and since \(c_kc'_l\gcd(k,l)\) is a non-negative integer for all \(k\) and \(l\), this proves the theorem.
\end{proof}

\begin{corollary}
	\label{lem:part2}
	Let \(D \in \GL_n(\ZZ)\) and \(d \in \ZZ^n\). The function
	\begin{equation*}
	\exp \sum_{k=1}^{\infty}
	O\left(\I_n-D^k,\left[\sum_{i=0}^{k-1}D^i\right] d\right) \frac{z^k}{k}
	\end{equation*}
	is a rational function.
\end{corollary}
\begin{proof}
From \cref{thm:Oissumsequences} we know that the \(O(...)\) equals a finite sum \(c_1a^1 + \cdots + c_la^l\). Hence:
\begin{align*}
	\exp \sum_{k=1}^{\infty}
O\left(\I_n-D^k,\left[\sum_{i=0}^{k-1}D^i\right] d\right) \frac{z^k}{k} &= \exp \sum_{k=1}^{\infty} \left[ \sum_{i=1}^l c_i a^i \right]_k\frac{z^k}{k}\\
&= \exp \sum_{k=1}^{\infty} \left[ \sum_{i=1}^l c_i a^i_k \right]\frac{z^k}{k}\\
&= \exp \sum_{i=1}^l c_i \left[\sum_{k=1}^{\infty}a^i_k \frac{z^k}{k}\right]\\
&= \prod_{i=1}^l  \exp \left[ - c_i \log(1-z^i) \right]\\
&= \prod_{i=1}^l  (1-z^i)^{-c_i},
\end{align*}
which is a rational function.
\end{proof}

%%%%%%%%%%%%%%%%%%%%%%%%%%%%%%%%%%
\iffalse
\begin{lemma}
	Let \(D \in \GL_n(\ZZ), d \in \ZZ^n\). The quantity
	\begin{equation*}
	O\left(\mathds{1}_n-D^k,\left[\sum_{i=0}^{k-1}D^i\right] d\right)
	\end{equation*}
	is periodic in \(k\).
\end{lemma}
\begin{proof}
	Let us consider everything mod \(2\). Because \(\GL_n(\ZZ_2)\) is finite, there exists some \(l \in \NN\) such that \(D^l = \I_n\). Then clearly \(\I_n - D^l\) is periodic in \(k\) with period \(l\). On the other hand, we have 
	\begin{align*}
	\sum_{i=0}^{2l-1}D^i &= \sum_{i=0}^{l-1}D^i + \sum_{i=l}^{2l-1}D^i\\
	&= \sum_{i=0}^{l-1}D^i + \sum_{i=0}^{l-1}D^{i+l}\\
	&= \sum_{i=0}^{l-1}D^i + \sum_{i=0}^{l-1}D^i\\
	&= 0,
	\end{align*}
	hence this sum is periodic with period $2l$. So clearly $O(...)$ is periodic with period $2l$. 
\end{proof}\fi

\subsection{The formula}
With both terms taken care of, we can now state the following theorem:
\begin{theorem}
	Let \(\varphi\) be an automorphism of the group \(\Gamma = \langle \ZZ^n, (0,-\I_n)\rangle\). Then there exist \(a,b,l \in \NN_0\), \(\mu_1, \dots, \mu_a, \nu_1, \dots, \nu_b \in \CC\) and \(c_1, \dots, c_l \in \NN_0\) such that
	\begin{equation*}
	R_z(\varphi) = \frac{ \prod_{i=1}^{b}(1-\nu_i z)}{\prod_{i=1}^{a}(1-\mu_i z) \prod_{i=1}^l (1-z^i)^{c_i}}.
	\end{equation*}
	The radius of convergence \(r\) of this function is given by 
	\begin{equation*}
	r = \frac{1}{\max\{1, |\mu_1|, \dots, |\mu_a|, |\nu_1|, \dots,|\nu_b|\}}.
	\end{equation*}
\end{theorem}

\section{Crystallographic groups with diagonal holonomy \(\ZZ_2\)}
We will now decompose a general crystallographic group \(\Gamma\) with diagonal holonomy \(\ZZ_2\) into the direct product of a group \(\Gamma_1\) from the family discussed in the previous section with some power of \(\ZZ\), which we will call \(\Gamma_2\). This will allow us to write the Reidemeister zeta function of \(\Gamma\) in terms of Reidemeister zeta functions of \(\Gamma_1\).
\begin{theorem}
	A Reidemeister zeta function of a crystallographic group with diagonal holonomy \(\ZZ_2\) is rational.
\end{theorem}
\begin{proof}
Any \(n\)-dimensional crystallographic group \(\Gamma\) with diagonal holonomy  \(\ZZ_2\) is of the form 
\begin{equation*}
\Gamma \cong \left\langle \ZZ^n, \left(a,\begin{pmatrix}
-\I_k & 0\\0 & \I_{n-k}\end{pmatrix}\right)\right\rangle.
\end{equation*}
If \(\Gamma\) is torsion-free, it is a Bieberbach group and hence its Reidemeister zeta functions (if any) are rational. On the other hand, if \(\Gamma\) is not torsion-free, then \(\Gamma\cong \ZZ^n \rtimes \ZZ_2\) and from this, it follows that we may assume without loss of generality that \(a = 0\) and hence
\begin{equation*}
 \Gamma \cong\underbrace{ \left\langle \ZZ^k, \left(0,-\I_k\right)\right\rangle}_{=\Gamma_1} \times \underbrace{\ZZ^{n-k}\vphantom{
 \left\langle \ZZ^k, \left(0,-\I_k\right)\right\rangle}}_{= \Gamma_2}.
\end{equation*}
In this decomposition \(\Gamma_2\) is the centre of \(\Gamma\), so is a characteristic subgroup of \(\Gamma\), but it is also not too hard to check that \(\Gamma_1\) is a characteristic subgroup as well.
If \(k = 0\) then \(\Gamma \cong \ZZ^n\) which is Bieberbach, and if \(k = 1\) then \(\Gamma / \langle e_2, e_3, \dots, e_n \rangle \cong \langle \ZZ, (0,-1)\rangle\), which has the \(R_\infty\)-property, hence \(\Gamma\) then has the \(R_\infty\)-property too. Of course, if \(n = k\), we simply have a group of the family we studied in \cref{sec:Zn0-1}. We thus assume that \(2 \leq k \leq n-1\). For any \(\varphi \in \Aut(\Gamma)\), we have that \(\varphi = \varphi_1 \times \varphi_2\) where \(\varphi_i\) ($i=1,2$) is an automorphism of $\Gamma_i$. We already know what \(R_{\varphi_1}(z)\) looks like, so we focus on \(\varphi_2\).

We know that \(\varphi_2 \in \Aut(\ZZ^{n-k}) \cong \GL_{n-k}(\ZZ)\), so let us denote its eigenvalues by \(\lambda_1, \dots, \lambda_{n-k}\). Let \(l \in \NN\) and assume that \(R(\varphi_2^l) < \infty\), then
\begin{equation*}
R(\varphi_2^l) = \left|\det\left(\I_{n-k}-\varphi_2^l\right)\right| = \left|\prod_{i=1}^{n-k} (1-\lambda_i^l) \right| = \prod_{i=1}^{n-k}|1-\lambda_i^l|.
\end{equation*}
Due to the same reasons as used in the proof of \cref{lem:firsttermkthpowers}, there exist \(a,b \in \NN_0\), \(\mu_1, \dots, \mu_a, \nu_1, \dots, \nu_b \in \CC\) such that
\begin{equation*}
R(\varphi_2^l) = \sum_{i=1}^a  \mu_i^l - \sum_{i=1}^b \nu_i^l.
\end{equation*}
Therefore, we have that
\begin{align*}
R_\varphi(z) &= \exp \sum_{k=1}^\infty R(\varphi^k) \frac{z^k}{k}\\
&= \exp \sum_{k=1}^\infty  R(\varphi_1^k)R(\varphi_2^k) \frac{z^k}{k}\\
&=  \exp \sum_{k=1}^\infty  R(\varphi_1^k) \left[\sum_{i=1}^a  \mu_i^k - \sum_{i=1}^b \nu_i^k\right]\frac{z^k}{k}\\
&= \exp  \left[\sum_{i=1}^a\sum_{k=1}^\infty R(\varphi_1^k) \frac{(\mu_iz)^k}{k}  -\sum_{i=1}^b \sum_{k=1}^\infty R(\varphi_1^k) \frac{(\nu_iz)^k}{k} \right]\\
&= \frac{\prod_{i=1}^a R_{\varphi_1}(\mu_iz)}{\prod_{i=1}^b R_{\varphi_1}(\nu_iz)},
\end{align*}
which is rational since \(R_{\varphi_1}\) is rational. If \(r\) is the radius of convergence of  \(R_{\varphi_1}\), then the radius of convergence of \(R_\varphi\) is given by \(r / \max \{|\mu_i|, |\nu_i|\}\).
\end{proof}

As a consequence of the above and the discussion of section 3 we may conclude with the following  theorem

\begin{theorem}
	A Reidemeister zeta function of an almost-crystallographic group of dimension at most \(3\) is rational.
\end{theorem}
%\bibliographystyle{plain}
%\bibliography{ReidemeisterSpectrum_Tertooy_Sam_AM_ArXiv_Corr}
\printbibliography
\end{document}